\documentclass[10pt,leqno]{amsart}
\usepackage{graphicx}
\usepackage{indentfirst,csquotes}
\usepackage[indentafter]{titlesec}
\titleformat{name=\section}{}{\thetitle.}{0.8em}{\centering\scshape}
\titleformat{name=\subsection}[runin]{}{\thetitle.}{0.5em}{\bfseries}[.]
\titleformat{name=\subsubsection}[runin]{}{\thetitle.}{0.5em}{\itshape}[.]
\titleformat{name=\paragraph,numberless}[runin]{}{}{0em}{}[.]
\titlespacing{\paragraph}{0em}{0em}{0.5em}
\titleformat{name=\subparagraph,numberless}[runin]{}{}{0em}{}[.]
\titlespacing{\subparagraph}{0em}{0em}{0.5em}
\usepackage{lipsum}

\usepackage{amssymb,amsthm,amsmath}
\usepackage{xcolor,paralist,hyperref,titlesec,fancyhdr,etoolbox,diagbox,multirow}
\usepackage[all]{xy}
\theoremstyle{definition}
\newtheorem{definition}{Definition}[]

\newtheorem{example}[definition]{Example}
\theoremstyle{plain}
\newtheorem{theorem}[definition]{Theorem}
\newtheorem{corollary}[definition]{Corollary}
\newtheorem{conjecture}[definition]{Conjecture}
\newtheorem{lemma}[definition]{Lemma}
\newtheorem{proposition}[definition]{Proposition}
\renewenvironment{proof}{\noindent\textsc{Proof.}\quad}{\qed}

\hypersetup{ colorlinks=true, linkcolor=black, filecolor=black, urlcolor=black }

\newcommand\Z{\mathop{}\!\mathbb{Z}}
\newcommand\Q{\mathop{}\!\mathbb{Q}}

\newcommand\p{\mathop{}\!\mathfrak{p}}
\newcommand\m{\mathop{}\!\mathfrak{m}}
\newcommand\OK{\mathop{}\!\mathcal{O}}

\newcommand\Hom{\mathop{}\!\mathrm{Hom}}
\newcommand\Spec{\mathop{}\!\mathrm{Spec}}
\newcommand\ch{\mathop{}\!\mathrm{ch}}
\newcommand\td{\mathop{}\!\mathrm{td}}

\newcommand\CH{\mathop{}\!\mathrm{CH}}

\newcommand\rank{\mathop{}\!\mathrm{rank}}
\newcommand\X{\mathop{}\!\mathcal{X}}
\newcommand\Y{\mathop{}\!\mathcal{Y}}
\newcommand\M{\mathop{}\!\mathcal{M}}
\begin{document}
\title{Riemann-Hurwitz Formula for Arithmetic Surfaces} 
\author[Subsun Jupp]{Ziyang ZHU}
\date{\today}
\address{School of Mathematical Sciences, Capital Normal University, Beijing 100048, China}
\email{zhuziyang@cnu.edu.cn}
\maketitle

\let\thefootnote\relax

\begin{abstract}
In this paper, we presents a method for factoring morphisms between arithmetic surfaces based on the regularity of arithmetic surfaces. Using this factorization, we derive a Riemann-Hurwitz formula satisfied by the ramification divisor and the canonical divisor on arithmetic surfaces. We also extend this formula to Arakelov theory.
\end{abstract}

\bigskip

\section{Introduction}\label{s1}

In complex algebraic geometry, surjective morphisms between smooth projective curves (Riemann surfaces) are all branched covering maps, and their branching behavior can be described by the following Riemann-Hurwitz formula
\[2\cdot\mathrm{genus}(X)-2=\deg(\varphi)(2\cdot\mathrm{genus}(Y)-2)+\sum_{x\in X}(e_x-1),\quad\varphi:X\to Y.\]
Here, the ramification index $e_x$ is actually the degree of the tangent map of $\varphi$ at $x$. In number theory, we also know that the integral extension of some ring of algebraic integers (arithmetic curve) also carries information about ramifications, with the corresponding Riemann-Hurwitz formula being
\[-\chi(\OK_L)=-[L:K]\cdot\chi(\OK_K)+\frac{1}{2}\deg(\mathrm{codifferent}(L/K)),\quad\OK_K\subseteq\OK_L,\]
where $\chi$ is the arithmetic Euler characteristic. In fact, through the adjunction formula, the theories in these two different frameworks can be well compatible. For example, any extension of $\mathbb{P}_{\mathbb{C}}^1$ or $\Z$ must be ramified, and one can use (profinite) Galois groups, or (\'{e}tale) fundamental groups, to study the relationship between ramification points and rational functions.

All the conclusions above pertain to one-dimensional curves. For higher-dimensional cases, there are not many results more detailed than the adjunction formula. The most direct difficulty is the lack of a good way to describe ramification divisors, especially when the base scheme is not a point and the relative dimension is greater than $0$. There has been some work, for example, Iversen \cite[Theorem IV.1]{ive} proved that for a map $\varphi:X\to Y$ with finite fibers between two non-singular algebraic surfaces defined over an algebraically closed field of characteristic zero, there is
\[\chi(X)=\deg(\varphi)\cdot\chi(Y)+\sum_ib_i(2p_g(B_i)-2)+\sum_{y\in Y}\lambda_{\varphi}(y),\]
where $\chi$ is Euler number (degree of the second Chern class) and $\sum_ib_iB_i$ is the branch divisor on $Y$. The term $\lambda_{\varphi}(y)$ is a local contribution which is zero outside the singular points of the branch locus. Additionally, Sun \cite{sun1, sun2} studied the corresponding formula where $\varphi$ is the structure morphism of an arithmetic surface over a complete discrete valuation ring. In his setting, the so-called ramification divisor of $\varphi$ is constructed by the Fitting ideal sheaf; the ramification number of $\varphi$ is understood as the Artin conductor of $\varphi$.

We also start from the adjunction formula to study the Riemann-Hurwitz formula on arithmetic surfaces. Roughly speaking, we have:
\subsection*{Theorem}
Let $\varphi:\X\to\Y$ be a finite flat separable morphism of arithmetic surfaces over $\Z$. Under some assumption, there exist canonical divisors $K_{\X/\Z}$ and $K_{\Y/\Z}$ such that the following equality of Weil divisors holds
\[K_{\X/\Z}=\varphi^*K_{\Y/\Z}+R_{\mathrm{geo}}+R_{\mathrm{ari}},\]
where $R_{\mathrm{geo}}$ is horizontal and $R_{\mathrm{ari}}$ is vertical.\\

The reason we express this only in terms of divisors is that there are no good analogs of degree and genus on arithmetic surfaces. We have thoroughly explored the ramification parts $R_{\mathrm{geo}}$ and $R_{\mathrm{ari}}$ in the above formula, which arise from the factorization of the morphism $\varphi$:
\subsection*{Theorem}
Assume $\varphi$ as above. Suppose $\X$ is a $\OK_K$-scheme and $\Y$ is a $\Z$-scheme. Then $\varphi$ can be factored through an arithmetic surface $\M:=\Y\times_{\Z}\OK_K$ by the universal property of pull-back
\[\varphi:\X\to\M\to\Y,\]
such that $\M\to\Y$ may only branch over those irreducible (entire) vertical fibers.\\

Even though our factorization can handle the vast majority of cases, we still have reason to believe that exceptional cases will not occur (Conjecture \autoref{20}). It is important to emphasize that addressing this technical aspect strongly relies on the properties of the arithmetic surface itself; in other words, the existence of $\varphi$ inherently comes with many requirements.

In addition, we apply our formula and the Grothendieck-Riemann-Roch theorem to compute the first (arithmetic) Chern class of the direct image of the structure sheaf in the (arithmetic) Chow group.

In \S\ref{s2}, we review some basic concepts in the theory of arithmetic surfaces and discuss the ramification behavior of ring extensions. In \S\ref{s3}, we present an important theorem (Theorem \autoref{5}) that is implied by the regularity of arithmetic surfaces. In \S\ref{s4} and \S\ref{s5}, we apply the adjunction formula to provide a line bundle version of the Riemann-Hurwitz formula on arithmetic surfaces and offer a method to transition from line bundles to divisors or cycles. Here, we propose a conjecture that serves as a complement to Theorem \autoref{5}. In the final section \S\ref{s6}, we will apply the formula (and its Arakelov version) to perform some calculations.

\subsection*{Acknowledgments}
The author would like to express gratitude to Professor Fei XU for encouraging me to explore this problem and for dedicating time to discuss it with me. I also extend my thanks to Mingqiang FENG for his assistance in proving Proposition \autoref{9}. Additionally, I appreciate Xiang LI for his insightful remark that inspired the final section of this paper.

\section{The Geometry of Arithmetic Surfaces}\label{s2}
An arithmetic surface $\X$ is a regular, integral scheme of dimension $2$ together with a projective, flat and of finite type morphism $\X\to\Spec(\OK_K)$, where $\OK_K$ is the ring of integers of a number field $K$. The generic fiber $X:=\X\times_{\OK_K}K$ is a smooth, geometrically connected and projective curve over $K$. So $\X$ is a regular model for $X$ over $\OK_K$.

A cycle of codimension $1$, or a Weil divisor on $\X$, is a sum of some codimension $1$ irreducible closed subschemes of $\X$. It is well-known that an arithmetic surface $\X\to\Spec(\OK_K)$ is fibered, and the fibers that lie on closed points of $\Spec(\OK_K)$ are the union of curves over the residue fields. Furthermore, these fibers are connected by Zariski's connectedness theorem (see \cite[Proposition III.3.3]{lan}). The prime Weil divisors on $\X$ can be divided into the following two cases:
\begin{itemize}
\item Horizontal divisor: Zariski closure of some closed point of the generic fiber $X$ which is finite and surjective to $\Spec(\OK_K)$.

\item Vertical divisor: Irreducible component of some closed fiber $\X\times_{\Z}(\OK_K/\p)$ where $\p$ is a prime ideal of $\OK_K$.
\end{itemize}

We review how to obtain divisors from line bundles. Since our arithmetic surface $\X\to\Spec(\OK_K)$ is integral, the line bundles can be identified with Cartier divisor classes. That is, there is an isomorphism
\[\mathrm{CaCl}(\X)\overset{\sim}{\longrightarrow}\mathrm{Pic}(\X),\quad D=\{(U_i,f_i)\}\longmapsto\OK_{\X}(D).\]
The inverse map is trivial: given a line bundle (locally free sheaf of rank $1$) $\mathcal{L}$, one can find an open covering $\{U_i\}$ of $\X$ such that $\mathcal{L}|_{U_i}$ is free and generated by some $f_i$ for each $i$. It is easy to check that $D=\{(U_i,f_i)\}$ is a Cartier divisor such that $\mathcal{L}\cong\OK_{\X}(D)$. Thanks for the regularity of $\X$, at this point, Weil divisors and Cartier divisors can be converted to each other. To be more precise, there is an isomorphism
\[\mathrm{CaCl}(\X)\overset{\sim}{\longrightarrow}\mathrm{WeCl}(\X),\quad D=\{(U_i,f_i)\}\longmapsto\sum_{Z}\mathrm{mult}_Z(D)[Z],\]
where $Z$ runs through all codimension $1$ normal irreducible closed subschemes of $\X$.

In order to study the branching property of morphisms between arithmetic surfaces, we need to provide a definition of ramifications, which coincides with concepts familiar to us in algebraic geometry. Noting that this is a local property, we will first review their behavior in the local setting.

Let $\phi:(A,\m_A)\to(B,\m_B)$ be a local homomorphism of local rings. Since we only concern codimension $1$ normal cycles on Noetherian regular schemes, we may assume $A$ and $B$ are both integrally closed regular local rings with Krull dimension $1$, hence are discrete valuation rings. Suppose $(\pi_A)=\m_A$ and $(\pi_B)=\m_B$ are uniformizers. We say $\phi$ has ramification index $e$ and residue index $f$, if $B\phi(\pi_A)=(\pi_B)^e$ and $B/(\pi_B)$ is a degree $f$ finite separable extension over $A/(\pi_A)$.

\begin{definition}\label{1}
Let $\varphi:\X\to\Y$ be a finite flat separable morphism of arithmetic surfaces over $\Z$. Let $W$ be a prime Weil divisor on $\Y$ and $Z$ is an irreducible component of $\varphi^*W$, i.e. $\varphi(Z)=W$. Assume $Z$ and $W$ are normal curves, we have a local homomorphism $\varphi_Z^{\sharp}:\OK_{\Y,W}\to\OK_{\X,Z}$ of discrete valuation rings. We say $\varphi$ has ramification index $e_Z$ and residue index $f_Z$ at $Z$, if $\varphi_Z^{\sharp}$ has ramification index $e_Z$ and residue index $f_Z$. Now, define the ramification divisor of $\varphi$ on $\X$ to be
\[R:=\sum_{W}\sum_{Z\text{ irreducible, normal and }\varphi(Z)=W}(e_Z-1)\cdot Z,\]
where $W$ runs through all codimension $1$ normal irreducible closed subschemes of $\Y$. Obviously, this is a finite sum.
\end{definition}

Since our local rings are Dedekind, the "$efg$" equality in algebraic number theory also holds. That is, $e_Zf_Z=\deg(\varphi_Z^{\sharp})$.

In some specific cases, the discriminant can delineate the support of the ramification divisor.

\begin{proposition}\label{2}
Let $A$ be a unique factorization domain, $A[\theta]$ is a single ring extension with $\theta$ a root of some monic polynomial $f\in A[x]$ such that $A[\theta]$ is regular. Consider the discriminate of $f$,
\[D(f):=\prod_{i\neq j}(\theta_i-\theta_j)^2\in A,\]
where $\theta_i$ are roots of $f$. Then the ramification divisor of $\Spec(A[\theta])\to\Spec(A)$ is generated by some prime divisors of $D(f)$.
\end{proposition}
\begin{proof}
We shall find those principal prime ideal $\p$ of $A$ such that the homomorphism $A_{\p}\to A_{\p}[x]/f$ is ramified. That is, there exists a maximal ideal $\m$ of $A_{\p}[x]/f$ with $\p(A_{\p}[x]/f)_{\m}=\m^e$ for some $e>1$. However, consider the residue field $\kappa_{\p}:=A_{\p}/\p A_{\p}$ of $\p$, the image $\overline{f}\in\kappa_{\p}[x]$ of $f\in A_{\p}[x]$ is a polynomial on a field, then $\overline{f}$ has multiple roots since $\kappa_{\p}[x]/\overline{f}\cong A_{\p}[x]/(\p,f)$ has nilpotent elements. Hence, $D(\overline{f})=0$, which means $D(f)\in\p$.
\end{proof}

\begin{example}\label{3}
Not every prime divisor of the discriminant necessarily contributes to ramifications; such exceptional cases primarily occur in vertical divisors. For example, in the ring extension $\Z[x]\subseteq\Z[x][y]/(y^2-x)$, the discriminate is $4x$, but $(2)$ is not ramified.
\end{example}

Sometimes we need the following condition for morphisms. Note that any ring extension in this setting can be written by a tower of single extensions, since the characteristic is zero.

\begin{definition}\label{4}
Let $\varphi$ be a finite morphism of arithmetic surfaces. Suppose $\varphi$ can be described locally by a ring extension $A\to B$ with a tower of single extensions
\[A\subseteq A[\theta_1]\subseteq\cdots\subseteq A[\theta_1,\cdots,\theta_{n-1}]\subseteq A[\theta_1,\cdots,\theta_{n-1},\theta_n]=B.\]
We say that $\varphi$ ramifies horizontally if the ramification divisors of every single extension above
\[\Spec(A[\theta_1])\to\Spec(A),~\cdots,~\Spec(A[\theta_1,\cdots,\theta_n])\to\Spec(A[\theta_1,\cdots,\theta_{n-1}])\]
are all horizontal divisors.
\end{definition}

Although this definition may seem complicated, it is not that significant. On one hand, it merely avoids the possibility of branching simultaneously over horizontal and vertical divisors; on the other hand, we conjecture that such situations are unlikely to occur (see Conjecture \autoref{20}).

\section{Factorization of Morphisms}\label{s3}
In this section, we study morphisms between arithmetic surfaces. We will find that once such a morphism can be expressed, it inadvertently imposes restrictions on the geometric or arithmetic properties of the arithmetic surfaces, and these restrictions can sometimes lead to strong conclusions. The main theorem of this section is as follows:

\begin{theorem}\label{5}
Let $\varphi:\X\to\Y$ be a finite flat separable morphism of arithmetic surfaces, where $\X$ is a $\OK_K$-scheme and $\Y$ is a $\Z$-scheme. Then $\varphi$ can be factored through an arithmetic surface $\M:=\Y\times_{\Z}\OK_K$ by the universal property of pull-back, such that $\beta$ may only branch over those irreducible vertical fibers.
\[\xymatrix{\X\ar[dr]\ar@{-->}[r]^{\alpha}\ar@/^2pc/[rr]^{\varphi}&\M\ar[d]\ar[r]^{\beta}&\Y\ar[d]\\&\OK_K\ar[r]&\Z}\]
\end{theorem}

The regularity of arithmetic surfaces can provide a lot of valuable information for proving our theorem. Since regularity is a local property, we can start with affine schemes over discrete valuation rings. Let us consider $U=\Spec(A)$ with $A=R[t_1,\cdots,t_n]/(f_1,\cdots,f_r)$, where $R$ is a discrete valuation ring with a uniformizer $\pi$. We always assume that $\pi\notin(f_1,\cdots,f_r)$.

First, we establish the Jacobian criterion of regularity of closed points of $U$. Let $\m_x=(g_1,\cdots,g_n,\pi)$ be a maximal ideal of $A$ with corresponding closed point $x$ of $U$. Since it is easy to verify that the matrix $\left(\frac{\partial g_i}{\partial t_j}(x)\right)$ is invertible as a matrix over the residue field of $x$, we can define the Jacobian matrix $\widetilde{J}_x$ as
\[\widetilde{J}_x:=\left(J_x\quad\bigg|\quad\frac{f_i(x)}{\pi}\right),\]
where $J_x=\left(\frac{\partial f_i}{\partial g_j}(x)\right)$ is the unique solution of the following equation of matrices:
\[\left(\frac{\partial f_i}{\partial t_j}(x)\right)=\left(\frac{\partial f_i}{\partial g_j}(x)\right)\left(\frac{\partial g_i}{\partial t_j}(x)\right).\]
This definition is indeed the implicit function theorem.

\begin{proposition}\label{6}
$U$ is regular at some closed point $x$ if and only if $\rank(\widetilde{J}_x)=n+1-\dim\OK_{U,x}$.
\end{proposition}
\begin{proof}
Recall that $U$ is regular at $x$ means $\dim_{\kappa_x}\m_x/\m_x^2=\dim\OK_{U,x}$, where $\kappa_x$ is the residue field of $x$. Suppose
\[f_i=\sum_jh_{ij}g_j+\widehat{f_i},\]
where $\widehat{f_i}\in R\pi$. Let $\overline{h_{ij}}$ denote the element corresponding to $h_{ij}$ in $\kappa_x$, so the dimension of $\m_x/\m_x^2$ is $n+1-\rank(\overline{h_{ij}},\widehat{f_i}(x)/\pi)$ as a $\kappa_x$-linear space. On the other hand, the differentiation rule provides
\[\frac{\partial f_i}{\partial t_j}=\sum_k\left(h_{ik}\frac{\partial g_k}{\partial t_j}+g_k\frac{\partial h_{ik}}{\partial t_j}\right),\]
modulo $\m_x$ on both sides we obtain
\[\frac{\partial f_i}{\partial t_j}(x)=\sum_k\overline{h_{ik}}\frac{\partial g_k}{\partial t_j}(x).\]
By definition we have $\overline{h_{ij}}=\frac{\partial f_i}{\partial g_j}(x)\in\kappa_x$, which means $J_x=(\overline{h_{ij}})$.
\end{proof}

We provide an example to illustrate this proposition.

\begin{example}\label{7}
Let $R=\Z[x]/(x^3+x+3)$ and $\m=(x^2+1,3)$. Denote $f=x^3+x+3$ and $g=x^2+1$, one can compute
\[\frac{\partial g}{\partial x}\equiv\big(2\sqrt{-1}\big)\quad\text{and}\quad\frac{\partial f}{\partial x}\equiv(-2)\quad\text{modulo }\m\]
as matrices. Hence, $\widetilde{J}_x=\left(\sqrt{-1},1\right)$, which has rank $1$ as a matrix over $\mathbb{F}_3[\sqrt{-1}]$. So $R$ is regular at $\m$.
\end{example}

Let $(R,\pi)$ be a discrete valuation ring with fractional field $K$, let $L/K$ be a finite separated extension and $(S,\pi')$ be the integral closure of $R$ in $L$. Let $U=\Spec(R[t_1,\cdots,t_n]/(f_1,\cdots,f_r))$ be an variety over $R$ as before, and $U_S:=U\times_RS$ be its base-change.

\begin{lemma}\label{8}
Let $x=(g_1,\cdots,g_n,\pi)$ be a closed point of $U$.
\begin{itemize}
\item If $L/K$ is unramified, then $x$ is regular implies every point of $U_S$ that lies on $x$ is regular.

\item If $L/K$ is ramified and $x$ is regular, then every point of $U_S$ that lies on $x$ is regular if and only if the system of linear equations
    \[\left(\frac{\partial f_i}{\partial g_j}(x)\right)\textbf{v}=\left(\frac{f_i(x)}{\pi}\right)\]
    has solutions, where $\textbf{v}$ is a column vector.
\end{itemize}
\end{lemma}
\begin{proof}
The first item is trivial. For the second one, let us assume $\pi=\pi'^e$ for some $e>1$, so $\widehat{f_i}\in R\pi$ implies $\frac{\widehat{f_i}}{\pi'}\in S\pi'$. Suppose $y$ is a point of $U_S$ that lies on $x$ with Jacobian matrix $\widetilde{J}_y$ and $\widetilde{J}_x$, where
\[J_x=\left(\frac{\partial f_i}{\partial g_j}(x)\right),~\widetilde{J}_x=\left(J_x\quad\bigg|\quad\frac{f_i(x)}{\pi}\right)\]
and
\[J_y=\left(\frac{\partial f_i}{\partial h_j}(x)\right),~\widetilde{J}_y=\left(J_y\quad\bigg|\quad\frac{f_i(x)}{\pi'}\right)=\left(J_y\quad\big|\quad0\right).\]
Here $(g_1,\cdots,g_n,\pi)\subseteq(h_1,\cdots,h_n,\pi')$. Since a point of a variety over a field and its base-change to some finite separated extension have the same regularity, so $\rank(J_x)=\rank(J_y)$. Therefore, $y$ is regular if and only if
\[\rank(\widetilde{J}_x)=\rank(\widetilde{J}_y)=\rank(J_y)=\rank(J_x),\]
which leads to the conclusion.
\end{proof}

We emphasize that in the setting of the second case of Lemma \autoref{8}, i.e. when $L/K$ is ramified, the regularity of a fiber is completely determined by the information on the original variety and is independent of $L$.

\begin{proposition}\label{9}
Let $(R,\pi)$ be a discrete valuation ring with fractional field $K$ and residue field $k$, let $L/K$ be a finite separated extension and $(S,\pi')$ be the integral closure of $R$ in $L$, whose residue field is $l$. Let $U$ be a regular variety over $R$ and $U_S:=U\times_RS$ (resp. $U_k:=U\times_Rk$) be its base-change to the integral extension (resp. to the residue field). We have:
\begin{itemize}
\item If $L/K$ is unramified, then $U_S$ is also regular.

\item If $L/K$ is ramified, then $U_S$ is regular if and only if $U_k$ is regular.
\end{itemize}
\end{proposition}
\begin{proof}
The unramified case is trivial by Lemma \autoref{8}, so we only need to verify the remaining case. Our $U_S$ is regular if and only if
\[\rank(\widetilde{J}_x)=\rank(J_x)=\rank(J_{\overline{x}})\]
for all $x$ of $U$ and its reduction $\overline{x}$. This is equivalent to $U_k$ is regular, since $\dim\OK_{U_k,\overline{x}}=\dim\OK_{U,x}-1$ (here we need $\pi$ is not generated by $f_1,\cdots,f_r$).
\end{proof}

Indeed, in the ramified case, if $U_k$ is regular, then $U$ and $U_S$ have the same regularity. This property is also independent of the choice of the ramified extension $R\subseteq S$.

\begin{example}\label{10}
Let $\Y=\Spec(\Z[x,y]/(xy-2))$ and consider $\Y_S=\X\times_{\Z}\Z[\sqrt{-1}]$, where $S=\Spec(\Z[\sqrt{-1}])$ is the spectrum of a Dedekind domain with discriminate $\Delta_S=-4$. Now $\Y_S$ is regular everywhere except on the fiber of $(2)\subseteq\Z$ by Lemma \autoref{8}. In fact, $\Y_S$ is not regular at $(1+\sqrt{-1},x,y)$, since after modulo $2$ the point $(x,y)$ of $\Y\times_{\Z}\mathbb{F}_2$ is the unique point that not regular.
\end{example}

The following corollary states that the ramification divisor on a base-change scheme of given arithmetic surface over $\Z$ by $\OK_K$ must be entire fibers.

\begin{corollary}\label{11}
Let $\Y$ be an arithmetic variety over $\Z$, and let $\OK_K$ be some ring of integers (with discriminate $\Delta_K$) such that $\Y\times_{\Z}\OK_K$ is regular. Then, the irreducible closed subscheme $\Y\times_{\Z}\mathbb{F}_p$ generates all the vertical divisors of $\X$ that lie on $p$ if $p|\Delta_K$.
\end{corollary}
\begin{proof}
By Proposition \autoref{9}, a necessary condition for $\Y\times_{\Z}\OK_K$ to be regular is that $\Y\times_{\Z}\mathbb{F}_p$ is irreducible at the ramified primes.
\end{proof}

This explains a geometric property of $\beta$ in Theorem \autoref{5}. To study $\alpha$, we need the following lemma. Recall that we can reduce to the case of single extensions, since the characteristic is zero.

\begin{lemma}\label{12}
Let $R$ be a discrete valuation ring and $A=R[t_1\cdots,t_n]/(f_1,\cdots,f_r)$. Consider a ring extension $A\subseteq A[t_{n+1}]/h$, suppose $x'$ is a closed point of $\Spec(A[t_{n+1}]/h)$ that lies on the close point $x$ of $\Spec(A)$, then
\[\rank(\widetilde{J}_{x'})\leq\rank(\widetilde{J}_x)+1.\]
The equality holds if $\frac{\partial h}{\partial t_{n+1}}(x')\neq0$.
\end{lemma}
\begin{proof}
Suppose $x'=(g_1',\cdots,g_n',g_{n+1}',\pi)$ and $x=(g_1,\cdots,g_n,\pi)$, their Jacobian matrices are $\widetilde{J}_{x'}$ and $\widetilde{J}_x$, respectively. Since we have
\[\left(\begin{array}{cc}\left(\frac{\partial f_i}{\partial t_k}(x')\right)_{r\times n}&(0)_{r\times1}\\\left(\frac{\partial h}{\partial t_k}(x')\right)_{1\times n}&\frac{\partial h}{\partial t_{n+1}}(x')\\\end{array}\right)=\left(\begin{array}{c}\left(\frac{\partial f_i}{\partial g_j'}(x')\right)_{r\times(n+1)}\\\left(\frac{\partial h}{\partial g_j'}(x')\right)_{1\times(n+1)}\\\end{array}\right)\left(\frac{\partial g_j'}{\partial t_k}(x')\right)_{(n+1)\times(n+1)},\]
and the matrix $\left(\frac{\partial g_j'}{\partial t_k}(x')\right)$ is invertible, so
\begin{align*}
\rank(\widetilde{J}_{x'})&=\rank\left[\left(\begin{array}{cc}\left(\frac{\partial f_i}{\partial g_j'}(x')\right)&\left(\frac{f_i(x')}{\pi}\right)\\\left(\frac{\partial h}{\partial g_j'}(x')\right)&\frac{h(x')}{\pi}\\\end{array}\right)\left(\begin{array}{cc}\left(\frac{\partial g_j'}{\partial t_k}(x')\right)&0\\0&1\\\end{array}\right)\right]\\
&=\rank\left[\left(\begin{array}{ccc}\left(\frac{\partial f_i}{\partial t_k}(x')\right)_{r\times n}&(0)_{r\times1}&\left(\frac{f_i(x')}{\pi}\right)_{r\times1}\\\left(\frac{\partial h}{\partial t_k}(x')\right)_{1\times n}&\frac{\partial h}{\partial t_{n+1}}(x')&\frac{h(x')}{\pi}\\\end{array}\right)\right].
\end{align*}
Note that $\left(\frac{\partial f_i}{\partial t_k}(x')\right)=\left(\frac{\partial f_i}{\partial t_k}(x)\right)=\left(\frac{\partial f_i}{\partial g_j}(x)\right)\left(\frac{\partial g_j}{\partial t_k}(x)\right)$, hence
\[\rank(\widetilde{J}_{x'})=\rank\left(\begin{array}{cc}\widetilde{J}_x&0\\\ast&\frac{\partial h}{\partial t_{n+1}}(x')\\\end{array}\right).\]
Now, the conclusion is directly provided by Gaussian elimination.
\end{proof}

According to this lemma, it can be obtained that the extension $\alpha$ can only break the regularity. That is to say, for $\alpha$ in Theorem \autoref{5}, if $\M$ is not regular, then $\X$ is also not regular by Proposition \autoref{6}.

\subsection*{Proof of Theorem \autoref{5}}
First, we show $\M$ is an arithmetic surface, and this only requires verifying the regularity of $\M$. Since regularity is a local property, we can observe it on affine charts. If $\M$ is not regular, by Lemma \autoref{12}, $\X$ is not regular, this is a contradiction. Now we have a regular $\OK_K$-scheme $\M=\Y\times_{\Z}\OK_K$, then by Proposition \autoref{9} and Corollary \autoref{11} we see that $\beta$ may only branch over those entire vertical fibers which lie on ramified primes of $\Z$.

\section{Ramification Theory: Geometric v.s. Arithmetic}\label{s4}
In original Riemann-Hurwitz formula on algebraic curves (see \cite[Proposition IV.2.3]{har} for a concrete example), the ramified information is recorded by the tangent bundle, so its determinant, the canonical line bundle.

\begin{definition}[{\cite[Definition 6.4.7]{liu}}]\label{13}
Let $\Y$ be a locally Noetherian scheme, and let $\varphi:\X\to\Y$ be a quasi-projective and local complete intersection morphism. Let $i:\X\to\mathcal{Z}$ be an immersion into a scheme $\mathcal{Z}$ that is smooth over $\Y$. We define the relative canonical line bundle of $\varphi:\X\to\Y$ to be the following line bundle on $\X$,
\[\omega_{\X/\Y}:=\det(\mathfrak{C}_{\X/\mathcal{Z}})^{\vee}\otimes i^*(\det\Omega_{\mathcal{Z}/\Y}^1),\]
where $\mathfrak{C}_{\X/\mathcal{Z}}$ is the conormal sheaf of $i$ and $\Omega_{\mathcal{Z}/\Y}^1$ is the sheaf of differentials of order $1$. This definition is independent of the choice of the decomposition $\X\to\mathcal{Z}\to\Y$ up to isomorphisms. Sometimes we write $\omega_{\X/\Y}$ as $\omega_{\varphi}$ to emphasize the morphism $\varphi$.
\end{definition}

Let $\varphi:\X\to\Spec(\OK_K)$ be an arithmetic surface, this $\varphi$ is automatically a local complete intersection morphism, so we have a relative canonical line bundle $\omega_{\varphi}$ on $\X$ by Definition \autoref{13}. We call any Cartier divisor $K_{\X/\Spec(\OK_K)}$ (or $K_{\varphi}$ for convenience) on $\X$ such that $\OK_{\X}(K_{\varphi})\cong\omega_{\varphi}$ a canonical divisor on $\X$ relative to $\Spec(\OK_K)$, where for any Cartier divisor $D=\{(U_i,f_i)\}$ the line bundle $\OK_{\X}(D)$ is the subsheaf of sheaf of stalks of rational functions on $\X$ such that $\OK_{\X}(D)|_{U_i}=f_i^{-1}\OK_{\X}|_{U_i}$ for all $i$. Such a divisor exists because $\X$ is integral (see \cite[Corollary 7.1.19]{liu}). We will subsequently present the specific construction of this divisor and a suitable cycle (Weil divisor) associated with it. For this reason, we need the following proposition \cite[Corollary 6.4.14]{liu}, which allows us to compute the rational sections of $\omega_{\varphi}$ explicitly.

\begin{proposition}\label{14}
Let $\Y=\Spec(A)$ be a Noetherian integral scheme, and let $\X$ over $\Y$ be an integral closed subscheme of $\Spec(A[T_1,\cdots,T_n])$ defined by an ideal generated by a regular sequence $F_1,\cdots,F_r$ with $r\leq n$. Suppose that
\[\Delta:=\det\left(\frac{\partial F_i}{\partial T_j}\right)_{1\leq i,j\leq r}\]
is non-zero in the function field $k(\X)$. Let $\eta$ be the generic point of $\X$.
\begin{itemize}
\item Let $t_i$ be the image of $T_i$ in $\OK_{\X}(\X)$, then $\omega_{\X/\Y,\eta}=(dt_{r+1}\wedge\cdots\wedge dt_n)\OK_{\X,\eta}$.

\item As a subsheaf of $\omega_{\X/\Y,\eta}$, we have $\omega_{\X/\Y}=\Delta^{-1}\cdot(dt_{r+1}\wedge\cdots\wedge dt_n)\OK_{\X}$.
\end{itemize}
\end{proposition}

Let $\varphi:\X\to\Y$ be a finite flat separable morphism of arithmetic surfaces over $\Z$, where we assume $\X$ is a $\OK_K$-scheme and $\Y$ is a $\Z$-scheme. By Theorem \autoref{5}, we have the following commutative diagram:
\[\xymatrix{\X\ar[drr]^{u}\ar@{-->}[rr]^{\alpha}\ar@/^2pc/[rrrr]^{\varphi}\ar@/_3pc/[drrrr]_{i}&&\Y\times_{\Z}\OK_K\ar[d]^{v}\ar[rr]^{\beta}&&\Y\ar[d]^{j}\\&&\OK_K\ar[rr]^{\gamma}&&\Z}\]

\begin{proposition}[Riemann-Hurwitz, bundle version]\label{15}
There is an isomorphism of line bundles
\[\omega_i\cong\varphi^*\omega_j\otimes\omega_{\alpha}\otimes u^*\omega_{\gamma}.\]
\end{proposition}
\begin{proof}
According to the adjunction formula and the base-change property for relative canonical line bundles \cite[Theorem 6.4.9]{liu}, we have
\begin{align*}
\omega_i&\cong\omega_u\otimes u^*\omega_{\gamma}\\
&\cong\omega_{\alpha}\otimes\alpha^*\omega_v\otimes\alpha^*v^*\omega_{\gamma}\\
&\cong\omega_{\alpha}\otimes\alpha^*(\omega_v\otimes v^*\omega_{\gamma})\\
&\cong\omega_{\alpha}\otimes\alpha^*(\omega_{\beta}\otimes\beta^*\omega_j)\\
&\cong\omega_{\alpha}\otimes\alpha^*\omega_{\beta}\otimes\alpha^*\beta^*\omega_j\\
&\cong\omega_{\alpha}\otimes\alpha^*v^*\omega_{\gamma}\otimes\varphi^*\omega_j\\
&\cong\omega_{\alpha}\otimes u^*\omega_{\gamma}\otimes\varphi^*\omega_j,
\end{align*}
as desired.
\end{proof}

Proposition \autoref{15} is an analogy to the Riemann-Hurwitz formula for curves. We call $\omega_{\alpha}\otimes u^*\omega_{\gamma}$ the ramification line bundle associated with $\varphi$, where $\omega_{\alpha}$ is the mixed data and $u^*\omega_{\gamma}$ the arithmetic data. As our terminology suggests, a ramified morphism between arithmetic surfaces is characterized jointly by its mixed data (if $\alpha$ ramifies horizontally, we then pass to generic fibers and the mixed data degenerates to the geometric data) and its arithmetic data (ramification of vertical divisors). In the next section, we will describe the contributions of the geometric data and arithmetic data via divisors.

\begin{example}\label{16}
Let $L/K$ be a field extension of number fields, it induces a canonical $\Z$-morphism $\varphi:\Spec(\OK_L)\to\Spec(\OK_K)$. By the adjunction formula we have
\[\omega_{\OK_L/\Z}\cong\omega_{\varphi}\otimes_{\OK_L}\varphi^*\omega_{\OK_K/\Z}\cong\omega_{\varphi}\otimes_{\OK_L}(\OK_L\otimes_{\OK_K}\omega_{\OK_K/\Z})\cong\omega_{\varphi}\otimes_{\OK_K}\omega_{\OK_K/\Z},\]
where $\omega_{\varphi}$ is in fact the codifferent ideal of $L/K$ (which encodes the ramifications) and $\omega_{\OK_K/\Z}=\Hom_{\Z}(\OK_K,\Z)$ is the canonical module of $K$. This well-known Riemann-Hurwitz formula for arithmetic curves is a degenerate case of Proposition \autoref{15} in algebraic number theory, see \cite[Proposition 3.3.11]{neu}.
\end{example}

\section{The Arithmetic Riemann-Hurwitz Formula}\label{s5}
We have obtained the bundle version of the Riemann-Hurwitz formula (Proposition \autoref{15}). In this section, we will establish the cycle form of the Riemann-Hurwitz formula, since given a line bundle we can produce a Weil divisor class. However, we cannot pick Weil divisors in this divisor class arbitrarily, because the divisor that describes the branching property is fixed, the rational equivalence may shift the divisor carrying this information elsewhere. Hence we need a canonical way to choose a suitable Weil divisor.

\begin{proposition}\label{17}
Under the assumptions of Theorem \autoref{5}, we assume that $\alpha$ ramifies horizontally. Let $R$ be the ramification divisor in Definition \autoref{1}, then $\OK_{\X}(R)$ is isomorphic to $\omega_{\alpha}\otimes u^*\omega_{\gamma}$. Moreover, $R$ can be decomposed as
\[R=R_{\mathrm{geo}}+R_{\mathrm{ari}},\]
where
\begin{itemize}
\item [(i).]$\OK_{\X}(R_{\mathrm{geo}})\cong\omega_{\alpha}$ and $\OK_{\X}(R_{\mathrm{ari}})\cong u^*\omega_{\gamma}$.

\item [(ii).]$R_{\mathrm{ari}}$ is vertical.

\item [(iii).]$R_{\mathrm{geo}}$ is horizontal.
\end{itemize}
\end{proposition}
\begin{proof}
Let $Y,M,X$ denote the generic fiber of $\Y,\M,\X$. We write $e(\square)_Z$ to emphasize the morphism $\square$, so $e(\varphi)_Z=e(\alpha)_Z\cdot e(\beta)_{\alpha(Z)}$. We have $e(\alpha)_Z=1$ for all $Z\cap Y=\varnothing$ since $\alpha$ ramifies horizontally; and $e(\beta)_{\alpha(Z)}=1$ for all $Z\cap Y\neq\varnothing$ by Theorem \autoref{5}. Define
\[R_{\mathrm{ari}}:=\sum_W\sum_{\varphi(Z)=W,Z\cap Y=\varnothing}(e(\varphi)_Z-1)\cdot Z=\sum_W\sum_{\beta\alpha(Z)=W,Z\cap Y=\varnothing}(e(\beta)_{\alpha(Z)}-1)\cdot Z\]
and
\[R_{\mathrm{geo}}:=\sum_W\sum_{\varphi(Z)=W,Z\cap Y\neq\varnothing}(e(\varphi)_Z-1)\cdot Z=\sum_W\sum_{\beta\alpha(Z)=W,Z\cap Y\neq\varnothing}(e(\alpha)_Z-1)\cdot Z,\]
which are two disjoint components of $R$. It is obvious that $R_{\mathrm{ari}}$ is vertical and $R_{\mathrm{geo}}$ is horizontal. Note that $\alpha$ is surjective, $Z\cap Y=\varnothing$ implies $\alpha(Z)\cap M=\varnothing$. So the summation indices for $R_{\mathrm{ari}}$ are those subschemes $Z$ of $\X$ such that $\alpha(Z)$ is vertical and ramified under $\beta$. But by Corollary \autoref{11} these $\alpha(Z)$ are entire vertical fibers that lie precisely on the ramified prime ideals under the ring extension $\Z\subseteq\OK_K$, and the ramification indices coincide with those in the arithmetic sense naturally. Hence, $\OK_{\X}(R_{\mathrm{ari}})\cong u^*\omega_{\gamma}$, where $\omega_{\gamma}=\Hom_{\Z}(\OK_K,\Z)$ is the $\OK_K$-module given by Example \autoref{16}.

Next, we show $\OK_{\X}(R_{\mathrm{geo}})\cong\omega_{\alpha}$. We claim that there is an one-to-one correspondence
\[\{\text{horizontal divisors on }\X\}\leftrightarrow\{\text{divisors on }X\},\quad D\mapsto i^*D,\]
where $i:X\to\X$ is the projection map in the pull-back diagram. Indeed, the inverse map is taking normalization of Zariski closure of given divisors on $X$. So $i^*R_{\mathrm{geo}}$ is the ramification divisor of $i^*(\alpha):X\to M$, which is fixed by $\omega_{i^*(\alpha)}\cong i^*\omega_{\alpha}$ by \cite[Theorem 6.4.9]{liu}. Thus, $R_{\mathrm{geo}}$ is rationally equivalent to a certain horizontal divisor determined by $\omega_{\alpha}$.
\end{proof}

According to Proposition \autoref{17}, we obtain the following corollary immediately.

\begin{corollary}[Riemann-Hurwitz, divisor version]\label{18}
Let $\varphi:\X\to\Y$ be a finite flat separable morphism of arithmetic surfaces over $\Z$. Under the assumptions of Proposition \autoref{17}, there exist canonical divisors $K_{\X/\Z}$ and $K_{\Y/\Z}$ such that the following equality of Weil divisors holds:
\[K_{\X/\Z}=\varphi^*K_{\Y/\Z}+R_{\mathrm{geo}}+R_{\mathrm{ari}}.\]
\end{corollary}

Here is an easy but computable example to present our results.

\begin{example}\label{19}
Let $\alpha:\X=\mathrm{Proj}(\Z[x,y,z]/(y^2-xz))\to\mathbb{P}_{\Z}^1$ be the natural $\Z$-morphism. First, we compute $\omega_{\X/\Z}$. By Proposition \autoref{14}, we observe that on the affine charts
\begin{center}
\begin{tabular}{c|c|c}
\hline
  $D_+(x)=\Spec\left(\frac{\Z[a,b]}{(b^2-a)}\right)$&$D_+(y)=\Spec\left(\frac{\Z[r,s]}{(1-rs)}\right)$&$D_+(z)=\Spec\left(\frac{\Z[u,v]}{(v^2-u)}\right)$\\
  $a=z/x,b=y/x$&$r=x/y,s=z/y$&$u=x/z,v=y/z$\\
  $\omega_{\X/\Z}(D_+(x))=\frac{da}{2b}\OK_{D_+(x)}$&$\omega_{\X/\Z}(D_+(y))=-\frac{dr}{r}\OK_{D_+(y)}$&$\omega_{\X/\Z}(D_+(z))=\frac{du}{2v}\OK_{D_+(z)}$\\
\hline
\end{tabular}
\end{center}
Hence, if on $D_+(x)$ we choose a section $\frac{1}{b}\frac{da}{2b}$, then it can be expressed as $-\frac{1}{v}\frac{du}{2v}$ (respectively, $-\frac{dr}{r}$) on the common part $D_+(x)\cap D_+(z)$ (respectively, $D_+(x)\cap D_+(y)$). This Cartier divisor defines a Weil divisor
\[K_{\X/\Z}=-[V_+(y,x)]-[V_+(y,z)].\]
Similarly, $\alpha^*K_{\mathbb{P}_{\Z}^1/\Z}=-2[V_+(y,x)]-2[V_+(y,z)]=2K_{\X/\Z}$. Therefore, $R_{\mathrm{geo}}=-K_{\X/\Z}=[V_+(y,x)]+[V_+(y,z)]$ is the ramification divisor of $\alpha$ by Corollary \autoref{18}. One can also compute this divisor by Proposition \autoref{14}.
\end{example}

\begin{conjecture}[Mixed Data $=$ Geometric Data]\label{20}
The condition "$\alpha$ ramifies horizontally" in Proposition \autoref{17} can be removed. That is, for any single ring extension which comes from a morphism of arithmetic surfaces, it cannot branch simultaneously over horizontal and vertical divisors.
\end{conjecture}

If the situation described in Conjecture \autoref{20} occurs, it is likely to disrupt the regularity of the arithmetic surface at the intersection of the horizontal and vertical divisors. A natural way to study this conjecture inevitably requires discussing the factorization of the discriminant ideal, which is quite challenging in general rings. It is of great significance to look for these cases or to describe the conditions under which this conjecture holds true.

\begin{example}\label{21}
Here is an inappropriate example. Consider the ring extension
\[\Z[x]\subseteq\Z[x][\sqrt{3x}+\sqrt{x+2}]=\Z[x][y]/f(x,y)\]
with minimal polynomial
\[f(x,y)=y^4+(-8x-4)y^2+(4x^2-8x+4).\]
The discriminate of $f$ is
\[D(f)=2^{14}\cdot3^2\cdot x^2\cdot(x-1)^2\cdot(x+2)^2.\]
Although we cannot discuss ramifications due to the breaking of regularity, we can still calculate the rank of Jacobian matrices of intersection points of divisors which appear in $D(f)$. One can see these ranks vanish everywhere.
\begin{center}
\begin{table}[!htbp]
\centering
\begin{tabular}{c|c|c|c|c}
\hline
\diagbox{fiber}{Jacobian}{$(x,y)$}&$\begin{array}{c}x=0\\(y^2-2)^2=0\end{array}$&$\begin{array}{c}x-1=0\\y^2=0\end{array}$&$\begin{array}{c}x-1=0\\y^2-12=0\end{array}$&$\begin{array}{c}x+2=0\\(y^2+6)^2=0\end{array}$\\
\hline
$(2)$&0&0&0&0\\
\hline
$(3)$&0&0&0&0\\
\hline
\end{tabular}
\end{table}
\end{center}
\end{example}

\section{The Grothendieck-Riemann-Roch Theorem}\label{s6}

The equality of divisors in Corollary \autoref{18} gives an equality of classes in $\CH^1(\X)$, the Chow group of $\X$, since $Z^1(\X)\twoheadrightarrow\CH^1(\X)$ is well-defined. That is, under the assumptions in \S\ref{s5}, the cycles
\[K_{\X/\Z}\sim\varphi^*K_{\Y/\Z}+R_{\mathrm{geo}}+R_{\mathrm{ari}}\]
are rational equivalence.

\begin{corollary}[Riemann-Hurwitz, weak]\label{22}
If in addition we assume $\X,\Y$ are smooth over $\Z$ and $\varphi$ is smooth and projective, then in $\CH^1(\X)_{\Q}:=\CH^1(\X)\otimes_{\Z}\Q$ we have
\[2c_1(\varphi_*\OK_{\X})\sim-\varphi_*(R_{\mathrm{geo}}+R_{\mathrm{ari}}).\]
\end{corollary}
\begin{proof}
By the Grothendieck-Riemann-Roch theorem (\cite{grs}, or see Proposition \autoref{24} below), we have
\[\varphi_*(\ch(\OK_{\X})\td(T_{\varphi}))=\ch(\varphi_*\OK_{\X}-R^1\varphi_*\OK_{\X}+\cdots),\]
where $T_{\varphi}$ is the relative tangent bundle of $\varphi$. The higher direct images $R^{>0}\varphi_*\OK_{\X}$ vanish by \cite[Proposition 5.2.34]{liu}. Note that $\ch=\rank+c_1+\cdots$ and $\td=1+\frac{c_1}{2}+\cdots$, the degree one part of the formula above is $\varphi_*(c_1(T_{\varphi}))=2c_1(\varphi_*\OK_{\X})$. Since $\varphi$ is smooth, there is an exact sequence
\[0\to\varphi^*\Omega_{\Y/\Z}\to\Omega_{\X/\Z}\to\Omega_{\X/\Y}\to0,\]
where $\Omega_{\X/\Y}=T_{\varphi}^{\vee}$ and $\Omega_{\X/\Z}=\omega_{\X/\Z},\Omega_{\Y/\Z}=\omega_{\Y/\Z}$ are locally free of rank one \cite[Corollary 6.2.6]{liu}. Hence, $c_1(\omega_{\X/\Z}^{\vee})=c_1(T_{\varphi})+c_1(\varphi^*\omega_{\Y/\Z}^{\vee})$, which means $c_1(T_{\varphi})=c_1(\varphi^*\omega_{\Y/\Z})-c_1(\omega_{\X/\Z})$. So we obtain
\[2c_1(\varphi_*\OK_{\X})=\varphi_*(c_1(\varphi^*\omega_{\Y/\Z})-c_1(\omega_{\X/\Z}))=\deg(\varphi)c_1(\omega_{\Y/\Z})-\varphi_*c_1(\omega_{\X/\Z}).\]
However, by Corollary \autoref{18} we have
\[\varphi_*K_{\X/\Z}=\deg(\varphi)K_{\Y/\Z}+\varphi_*(R_{\mathrm{geo}}+R_{\mathrm{ari}}).\]
By comparing the two equalities, the formula $2c_1(\varphi_*\OK_{\X})\sim-\varphi_*(R_{\mathrm{geo}}+R_{\mathrm{ari}})$ holds.
\end{proof}

This is an application of our Riemann-Hurwitz formula. We can also generalize this equality to Arakelov's arithmetic cycles \cite{ara} or \cite[Chapter III]{sabk}.

Let $\widehat{\CH}^1(\X)$ be the arithmetic Chow group of $\X$, defined as the quotient group $\widehat{Z}^1(\X)/\widehat{R}^1(\X)$, where
\[\widehat{Z}^1(\X):=\big\{(Z,g_Z):Z\in Z^1(\X),~g_Z\text{ a Green current for }Z(\mathbb{C})\big\}\]
with addition defined componentwise and $\widehat{R}^1(X)\subseteq\widehat{Z}^1(X)$ is the subgroup generated by pairs
\[\big(\mathrm{div}(f),-[\log|f_{\mathbb{C}}|^2]\big),\quad f\in k(\X)^{\times}.\]
The elements in $\widehat{\CH}^1(\X)$ are called Arakelov divisors (classes). By definition there is a natural map $Z^1(\X)\hookrightarrow\widehat{Z}^1(\X)\twoheadrightarrow\widehat{\CH}^1(\X)$, $D\mapsto[\widehat{D}]$. Based on this, we immediately obtain:

\begin{corollary}[Riemann-Hurwitz, Arakelov version]\label{23}
Under the assumptions of Corollary \autoref{18},
\[[\widehat{K_{\X/\Z}}]=\varphi^*[\widehat{K_{\Y/\Z}}]+[\widehat{R_{\mathrm{geo}}}]+R_{\mathrm{ari}}.\]
\end{corollary}

To obtain a conclusion similar to Corollary \autoref{22} in the Arakelov setting, we need the arithmetic Grothendieck-Riemann-Roch theorem.

\begin{proposition}[\cite{grs}]\label{24}
Let $\varphi:\X\to\Y$ be a flat projective $\Z$-morphism of arithmetic varieties, which is smooth over $\Q$. Let $K_0$ denote the Grothendieck group of vector bundles and $\widehat{K_0}$ denote the arithmetic Grothendieck group \cite[Definition 2.1]{grs}. Given a K\"{a}hler metric $g$ on $\X(\mathbb{C})$ which is invariant under the complex conjugation $F_{\infty}:\X(\mathbb{C})\to\X(\mathbb{C})$, then we have the following commutative diagram:
\[\xymatrix{
K_0(\X)\ar[dd]_{\sum(-1)^iR^i\varphi_*(\cdot)}\ar[rrrrr]^{\ch(\cdot)\td(T_{\varphi})}&&&&&\CH^{\cdot}(\X)_{\Q}\ar[dd]^{\varphi_*}&\\
&\widehat{K_0}(\X)\ar@{->>}[ul]\ar[dd]_{\widehat{\varphi}_*}\ar[rrrrr]^{\widehat{\ch}(\cdot)\widehat{\td}(\overline{T_{\varphi}})(1-\widehat{\ch}(\mathcal{R}(\overline{T_{\varphi,\mathbb{C}}})))}&&&&&\widehat{\CH}^{\cdot}(\X)_{\Q}\ar@{->>}[ul]\ar[dd]^{\varphi_*}\\
K_0(\Y)\ar[rrrrr]^{\ch(\cdot)}&&&&&\CH^{\cdot}(\Y)_{\Q}&\\
&\widehat{K_0}(\Y)\ar@{->>}[ul]\ar[rrrrr]^{\widehat{\ch}(\cdot)}&&&&&\widehat{\CH}^{\cdot}(\Y)_{\Q}\ar@{->>}[ul]}\]
where
\begin{itemize}
\item $\widehat{\ch},\widehat{\td}$ are arithmetic characteristic classes.

\item $\overline{T_{\varphi}}$ is the metrized relative tangent bundle on $\X$ whose metric induced by $g$.

\item $\mathcal{R}$ is the R-genus, where $\mathcal{R}(x):=\sum_{\text{odd }k\geq1}\left(\sum_{j=1}^k\frac{1}{j}+2\frac{\zeta'(-k)}{\zeta(-k)}\right)\zeta(-k)\frac{x^k}{k!}$.

\item $\widehat{\varphi}_*$ is the push-forward in arithmetic K-theory.
\end{itemize}
\end{proposition}

From the arithmetic Riemann-Roch theorem, we can obtain a weak Arakelov version of the Riemann-Hurwitz formula. We refer to this formula as weak because it is applicable in $\widehat{\CH}^1(\X)_{\Q}$ and, in most cases, it cannot be computed.

\begin{corollary}[Riemann-Hurwitz, Arakelov version, weak]\label{25}
Under the assumptions of Corollary \autoref{22}, we have
\[2\widehat{c_1}(\overline{\varphi_*\OK_{\X}})\sim\left(\varphi_*\widehat{\ch}(\mathcal{R}(\overline{T_{\varphi,\mathbb{C}}}))-1\right)\cdot\varphi_*\left([\widehat{R_{\mathrm{geo}}}]+R_{\mathrm{ari}}-\widehat{c_1}(\widetilde{\ch}(\mathfrak{E}_{\mathbb{C}}))\right)-2\widehat{c_1}(\varphi_{\overline{\OK_{\X}}}),\]
where
\begin{itemize}
\item The metric of $\overline{\varphi_*\OK_{\X}}$ is the Quillen metric (induced by the metric $g$ and the analytic torsion), see \cite[Definition 1.10]{bl}.

\item $\mathfrak{E}_{\mathbb{C}}$ is an exact sequence of metrized vector bundles on $\X(\mathbb{C})$, which is the pull-back of a sequence of metrized vector bundles
    \[\mathfrak{E}:\quad0\to\overline{T_{\varphi}}\to\overline{\Omega_{\X/\Z}}^{\vee}\to\varphi^*\overline{\Omega_{\Y/\Z}}^{\vee}\to0\]
    on $\X$ with metrics induced by $g$.

\item $\widetilde{\ch}$ is the Bott-Chern secondary characteristic form \cite[Theorem 1.29]{bgs}, which is defined on short exact sequences of metrized vector bundles.

\item $\displaystyle\varphi_{\overline{\OK_{\X}}}\in\bigoplus_{p\geq0}\frac{\{\omega\text{ smooth real $(p,p)$-form on $\X(\mathbb{C})$ s.t. }F_{\infty}^*\omega=(-1)^p\omega\}}{\mathrm{im}(\partial)+\mathrm{im}(\overline{\partial})}$ is the higher torsion form \cite{bk}, which satisfies the following equation
    \[\frac{i}{2\pi}\partial\overline{\partial}(\varphi_{\overline{\OK_{\X}}})=\int_{\X(\mathbb{C})/\Y(\mathbb{C})}\td(\overline{T_{\varphi,\mathbb{C}}})-\ch(\overline{\varphi_*\OK_{\X,\mathbb{C}}}).\]
\end{itemize}
\end{corollary}
\begin{proof}
By Proposition \autoref{24} we have $\varphi_*\big(\widehat{\td}(\overline{T_{\varphi}})(1-\widehat{\ch}(\mathcal{R}(\overline{T_{\varphi,\mathbb{C}}})))\big)=\widehat{\ch}(\widehat{\varphi}_*\overline{\OK_{\X}})$.
According to the construction of arithmetic Grothendieck group,
\[\widehat{c_1}(\overline{T_{\varphi}})+\widehat{c_1}(\varphi^*\overline{\omega_{\Y/\Z}}^{\vee})-\widehat{c_1}(\overline{\omega_{\X/\Z}}^{\vee})=\widehat{c_1}(\widetilde{\ch}(\mathfrak{E}_{\mathbb{C}})).\]
Hence, by the definition \cite[Section 3]{grs} of $\widehat{\varphi}_*$, the degree one part is
\[\varphi_*\Big(\big(\widehat{c_1}(\varphi^*\overline{\omega_{\Y/\Z}})-\widehat{c_1}(\overline{\omega_{\X/\Z}})+\widehat{c_1}(\widetilde{\ch}(\mathfrak{E}_{\mathbb{C}}))\big)(1-\widehat{\ch}(\mathcal{R}(\overline{T_{\varphi,\mathbb{C}}})))\Big)=2\big(\widehat{c_1}(\overline{\varphi_*\OK_{\X}})+\widehat{c_1}(\varphi_{\overline{\OK_{\X}}})\big),\]
since $\OK_{\X}$ is $\varphi$-acyclic. However,
\[\varphi_*\big(\widehat{c_1}(\varphi^*\overline{\omega_{\Y/\Z}})-\widehat{c_1}(\overline{\omega_{\X/\Z}})\big)\sim\varphi_*\varphi^*[\widehat{K_{\Y/\Z}}]-\varphi_*[\widehat{K_{\X/\Z}}]=-\varphi_*([\widehat{R_{\mathrm{geo}}}]+R_{\mathrm{ari}}).\]
By substituting and rearranging, we can obtain the final formula.
\end{proof}

\end{document}